\documentclass[11pt,reqno]{amsart}
\pdfoutput=1
 
\usepackage[table]{xcolor}
\usepackage[bmargin=2.5cm,tmargin = 3cm, hmargin=2cm]{geometry}
\usepackage{caption}

\usepackage[backref = page]{hyperref}
\hypersetup{
    colorlinks=true,
    linkcolor=blue,
    urlcolor =blue, 
    citecolor=magenta,
}

\usepackage{varwidth}

\renewcommand*{\backrefalt}[4]{%
\ifcase #1 %
\color{red} No citations.%
\or
(p.~#2).%
\else
(pp.~#2).%
\fi
}

\usepackage[nameinlink, capitalize, noabbrev]{cleveref}

\usepackage{amssymb,amsthm,amsmath, stmaryrd}

\usepackage{mathtools}
\usepackage{multirow}
\usepackage{float}
\usepackage{cancel}

\usepackage{booktabs}

\theoremstyle{plain}

\newtheorem{theorem}{Theorem}[section]

\crefname{thma}{Theorem}{Theorems}

\numberwithin{equation}{section}

\newtheorem{lemma}[theorem]{Lemma}

\theoremstyle{definition}
\newtheorem{remark}[theorem]{Remark}

\newcommand{\Z}{\mathbb{Z}}

\newcommand{\R}{\mathbb{R}}
\newcommand{\N}{\mathbb{N}}
\newcommand{\C}{\mathbb{C}}

\DeclareMathOperator{\SL}{SL}

\renewcommand{\Re}{\operatorname{Re}}

\renewcommand{\pmod}[1]{\mkern4mu(\mathrm{mod} \, #1)}

\makeatletter
\def\section{\@startsection{section}{1}%
  \z@{2.5\linespacing\@plus.5\linespacing\@minus.5\linespacing}{.5\linespacing}%
  {\normalfont\scshape\centering}}
\makeatother

\begin{document}

\title{On the asymptotics of certain colored partitions}

\author{Lukas Mauth}
\address{Department of Mathematics and Computer Science, Division of Mathematics, University of Cologne, \linebreak Weyertal 86-90, 50931 Cologne, Germany}
\email{lmauth@uni-koeln.de}

\keywords{partitions, modular forms, Dedekind eta function, Circle Method}

\subjclass[2020]{11B57, 11F03, 11F20, 11F30, 11F37, 11P82}

\begin{abstract}
We will prove an infinite family of asymptotic formulas for the logarithm of certain two-colored partitions. An infinite sub-family of these asymptotics was posed as a conjecture by Guadalupe.
\end{abstract}

\maketitle

\section{Introduction}
Let $n$ be a non-negative integer. A \emph{partition} of $n$ is a non-increasing sequence of positive integers that sum to $n$. The numbers which appear in this sequence are called the \emph{parts} of the partition. We define $p(n)$ to be the number of partitions of $n$, where we agree on the convention $p(0)=1$. Partitions and the variations thereof have been the subject of extensive research. One of the early milestones was the work of Hardy and Ramanujan \cite{Hardy-Ramanujan} in which they established the asymptotic formula
\[
p(n) \sim \frac{1}{4n\sqrt{3}}e^{\pi\sqrt{\frac{2n}{3}}},
\]
as $n \rightarrow \infty$, using their celebrated Circle Method. A crucial step in their method is to look at all $p(n)$ together in terms of their generating function and realize that it is essentially a modular form of weight $-\frac{1}{2}$. Indeed, the generating function satisfies the following relation \cite{Andrews}
\[
    P(q) \coloneqq \sum_{n=0}^{\infty} p(n)q^n = \prod_{n=1}^{\infty} \frac{1}{1-q^n},
\]
where the denominator of the right hand side is up to a $q$-power the Dedekind $\eta$-function, which is well-known to be a modular form of weigth $\frac{1}{2}.$
Rademacher improved the Circle Method and was able to find a formula for $p(n)$ in terms of a rapidly absolutely converging series involving only elementary functions \cite{Rademacher}. A short time after Rademacher's work, Zuckerman generalized \cite{Zuckerman} this result to state exact formulas of Rademacher-type (these exact formulas do not anymore necessarily involve only elementary functions) for the Fourier coefficients of any weakly holomorphic modular form of negative weight for any finite index subgroup of $\SL_2(\Z)$. This immediately gives a convenient way to find similar formulas for cousins of the partition function. For instance, let us define (see ) for any positive integer $r$ the function $a_r(n)$ by the relation
\[
    G_r(q) \coloneqq \sum_{n=0}^{\infty}a_p(n)q^n = \prod_{n=1}^{\infty} \frac{1}{(1-q^n)(1-q^{rn})} = P(q)P(q^r).
\]
This family of generating functions (complete or for specific values of $r$) has been considered in for instance \cite{AhmedBaruahDastidar, Guadalupe, Banerjee-Paule-Radu-Zeng}. Then, it is not hard too see that $a_r(n)$ counts the number of certain two-colored partitions of $n$, where a part is only allowed to have two colors if it is divisible by $r$ (compare \cite{AhmedBaruahDastidar}). In particular, the function $a_1(n)$ counts the number of two-colored partitions of $n$ without any conditions on the parts. The terminology of colored partitions (also known as multipartitions) together with many interesting properties of those partitions can for instance be found in \cite{Andrews2, Keith}. In \cite{Banerjee-Paule-Radu-Zeng} Banerjee, Radu, Paule and Zeng conjectured an asymptotic formula for $\log \left(a_2(n)\right)$ which the author established by giving an exact formula of Rademacher-type for $a_2(n)$, using Zuckerman's general result \cite{Mauth}. Guadalupe \cite{Guadalupe} generalized their conjecture to the claim that for all primes $p$ we have as $n \rightarrow \infty$ the following asymptotics
\[
\log \left(a_p(n)\right) \sim \pi \sqrt{\frac{2n(1+p^{-1})}{3}}-\frac{5}{4}\log \left(n\right) + \log \left(\frac{2\sqrt{3p}(1+p^{-1})^{\frac{3}{4}}}{24^{\frac{5}{4}}}\right) - \frac{c_p}{24\sqrt{6n}},
\]
where for any positive integer the constant $c_r$ is defined by
\[
    c_r \coloneqq \frac{135}{\pi\sqrt{1+r^{-1}}} + \pi \sqrt{\frac{(1+r)^3}{r}}.
\]
For $p=2$ this reduces to the aforementioned conjecture from \cite{Banerjee-Paule-Radu-Zeng}. Guadalupe has already shown \cite{Guadalupe} his conjecture for $p< 24$. He followed a similar approach the author used to prove the $p=2$ case in \cite{Mauth}. Indeed, he builds up upon the work of Sussman \cite{Sussman} who specialized Rademachers techniques for the Circle Method to a fairly large class of eta-quotients and obtained exact formulas for their Fourier coefficients as a convergent series of Rademacher-type. It is worthwile to mention that Chern \cite{Chern} established a very similar result for another class of eta-quotients. It turns out that for $p<24$ the function $G_p(q)$ satisfies the conditions of Sussman's result and this gave rise to a rapid computation of exact formulas for $a_p(n)$ from which an asymptotic formula for $\log \left(a_p(n)\right)$ easily follows. 

In this paper we prove the following Theorem.
\begin{theorem}\label{thm:main-theorem}
Let $r$ be a positive integer. Then, we have
\[
\log \left(a_r(n)\right) \sim \pi \sqrt{\frac{2n(1+r^{-1})}{3}}-\frac{5}{4}\log\left(n\right) + \log \left(\frac{2\sqrt{3r}(1+r^{-1})^{\frac{3}{4}}}{24^{\frac{5}{4}}}\right) - \frac{c_r}{24\sqrt{6n}}.
\]
\end{theorem}

\begin{remark}
We point out that \cref{thm:main-theorem} implies Guadalupe's conjecture \cite{Guadalupe} in full generality. We further point out that for $r=1$ this theorem recovers the asymptotics for the unrestricted two-colored partition function, see \cite[Corollary 2.2]{Iskander-Jain-Talvola}.
\end{remark}

Even though one could follow for the proof the same approach as for $p<24$ and provide exact formulas for $a_r(n)$, which is definitely feasible as the result of Zuckerman \cite{Zuckerman} shows, we will instead prove the asymptotics using a short Circle Method argument. We use the classical Hardy--Ramanujan Circle Method and compute the first few terms in the asymptotic expansion of $a_r(n)$. This is already enough to prove the conjecture in full generality. For simplicity we will first show the case when $r=p$ is a prime and afterwards establish the claimed asymptotics for any positive $r$.

\subsection*{Acknowledgements}
The author wishes to thank Badri Vishal Pandey for suggesting this project. Furthermore, the author thanks Kathrin Bringmann and Bernhard Heim for helpful discussions.

\section{Proof of Main Theorem}

Let $p$ be a prime. We will introduce the following setup for the Circle Method as described in \cite{Andrews}. Let $N= \sqrt{n}$ and assume throughout that $n > \frac{p+1}{24}$. We denote by $F_N$ the sequence of \emph{Farey fractions} up to order $N$. We define for adjacent Farey Fractions $\frac{h_1}{k_1} < \frac{h}{k} < \frac{h_2}{k_2}$ in $F_N$ the quantities\footnote{Note the typographical error in the definition of $\vartheta_{h,k}'$ and $\vartheta_{h,k}''$ in \cite{Andrews}. We use an equivalent representation of $\vartheta_{h,k}'$ and $\vartheta_{h,k}''$ (see \cite{RademacherJFunction}).}
\[
    \vartheta_{h,k}' \coloneqq \frac{1}{k(k_1+k)}, \quad \vartheta_{h,k}'' \coloneqq \frac{1}{k(k_2+k)}, \quad \vartheta_{0,1}=\frac{1}{N+1}.
\]
It is not hard to see that $\frac{1}{k_j+k}\leq \frac{1}{N}$ for $j=1,2$. Furthermore, we set for $-\vartheta_{h,k}' \leq \Phi \leq \vartheta_{h,k}''$ the variable $z = \frac{k}{n} - ik\Phi.$ Note that by construction we have $\Re(z) > 0$ as well as the bounds $\Re(1/z) \geq \frac{k}{2}$ and $|z| \ll \frac{1}{\sqrt{n}}$. We now write using Cauchy's integral formula
\[
a_p(n) = \frac{1}{2\pi i} \int_{C} \frac{G_p(n)}{q^{n+1}} dq,
\]
where $C$ is the circle of radius $r=e^{-\frac{2\pi}{N^2}}$ traversed once in the counterclockwise direction. The circle shall be parametrized by $q = e^{-\frac{2\pi}{N^2}+ 2 \pi i t}$, so the integral transforms into
\[
a_p(n) = \int_{0}^{1} G_p\left(e^{-\frac{2\pi}{N^2} + 2\pi i t}\right) e^{\frac{2\pi n}{N^2}-2\pi int} dt.
\]
We will now split the integral along the Farey arcs defined above and henceforth obtain
\begin{equation}\label{eqn:integral-farey-arcs}
a_p(n) = \sum_{\substack{0 \leq h < k \leq N \\ \gcd(h,k)=1}} e^{-\frac{2\pi ihn}{k}} \int_{-\vartheta_{h,k}'}^{\vartheta_{h,k}''}G_p\left(e^{\frac{2\pi i}{k}\left(h+iz\right)}\right) \cdot e^{\frac{2\pi n z}{k}} d\Phi.
\end{equation}

It is at this point where we will use the modular transformation law for $G_p(q)$. To state the transformation law we introduce for $\frac{h}{k} \in F_N$ and $\Re(z) > 0$ the variables $q=e^{\frac{2\pi i}{k}\left(h + iz\right)}$ and $q_1 = e^{\frac{2\pi i}{k}\left(h' + \frac{i}{z}\right)}$, where $h'$ is any solution to the congruence $hh' \equiv -1 \pmod{k}$. Recall that $G_p(q) = P(q)P(q^p)$. The modular transformation law for $P(q)$ is given by
\[
P(q) = \omega_{h,k}e^{\frac{\pi}{12k}(z^{-1}-z)}P(q_1),
\]
where $\omega_{h,k}$ is a certain $24k$-th root of unity, whose definition should not bother us (see \cite[Section 5]{Andrews} for more details). The only fact concerning $\omega_{h,k}$ we care about is that $\omega_{0,1} = 1$. There exist transformation laws for $P(q^p)$ (see for instance \cite{Bringmann-Mahlburg}). However, one will have to distinguish the two cases $p \mid k$ or $p \nmid k$. In the case where $p \mid k$ we have
\[
    P(q^p) = \omega_{h,\frac{k}{p}}\sqrt{z}e^{\frac{\pi p}{12k}(z^{-1}-z)}P(\zeta q_1^p)
\]
and in the other case we find
\[
    P(q^p) = \omega_{hp,k} \sqrt{pz} e^{\frac{\pi}{12k}((pz)^{-1}-pz)}P\left(\zeta q_1^{\frac{1}{p}}\right),
\]
where $\zeta = \zeta(p,h,k)$ is a certain root of unity whose exact value is not relevant for our computations.

After applying the transformation laws we will now split the sum in \eqref{eqn:integral-farey-arcs} into two parts, depending on whether $p \mid k$ or $p \nmid k$. We will start with the case that $p \nmid k$ since this will give us the dominating contribution from the sum, as we shall see. The sum in question is
\begin{equation}\label{eqn:contribution-p-coprime-to-k}
    \sum_{\substack{0 \leq h < k \leq N \\ \gcd(h,k)=1 \\
    p \nmid k}} e^{-\frac{2\pi ihn}{k} }\sqrt{p} \omega_{h,k}\omega_{hp,k} \int_{-\vartheta_{h,k}'}^{\vartheta_{h,k}''} z e^{\frac{\pi}{12k}\left(\left(1+\frac{1}{p}\right)\frac{1}{z}-(1+p)z\right)}e^{\frac{2\pi nz}{k}}P(q_1)P\left(\zeta q_1^{\frac{1}{p}}\right) d\Phi.
\end{equation}
At this point it will be convenient to plug in the $q$-expansions for $P(q_1)$. The integral becomes
\[
    \int_{-\vartheta_{h,k}'}^{\vartheta_{h,k}''} z \sum_{m_1,m_2 \geq 0} \zeta^{m_2} p(m_1)p(m_2) e^{\frac{2\pi ih'}{k}\left(m_1+p^{-1}m_2\right)} e^{\frac{\pi}{12k}(A(p,m_1,m_2)z^{-1}+(24n-p-1)z)} d\Phi,
\]
where we defined 
\[
A\coloneqq A(p,m_1,m_2) \coloneqq \left(1+\frac{1}{p}-24m_1-\frac{24m_2}{p}\right).
\]

Note that depending on the sign of $A$ the corresponding summand in the integral either grows or decays exponentially as $z \rightarrow 0$ (with the exception of $A=0$. However an easily calculation shows that this case can never occur). Therefore, we can expect that only those pairs $(m_1,m_2)$ will contribute where $A>0.$ First, we note that $m_1 \geq 1$ already implies $A(p,m_1,m_2)<0$. We denote $m_p$ as the minimal positive integer, such that $A(p,0,m_2) < 0.$

\begin{remark}
    	It is worthwhile to note that for $p<24$ one has $m_p = 0$. This is the exact point where the general result from \cite{Sussman} which was used in \cite{Guadalupe} to prove the conjecture for $p<24$ breaks. In fact Sussman himself mentiones in \cite{Sussman} that one of the conditions from his main result ensures that the principal part at each cusp of the modular form contains at most one summand. He mentiones as well that his proof can be easily adepted to deal with a more general situation. This goes in the direction of \cite{Zuckerman} who showed the most general form. The art here is to find an equilibrium between generality and applicability. In comparison between \cite{Mauth} and \cite{Guadalupe} the specialed result for eta-quotients yields much simpler an exact formula of Rademacher-type than Zuckerman's result.
\end{remark}

Indeed, we are now going to prove that the terms with $A<0$ do only have a polynomial contribution in $n$. We estimate the above integral using $\Re(1/z) \geq \frac{k}{2}$, $\Re(z) = \frac{k}{n}$, and $|z| \ll \frac{1}{\sqrt{n}}$ to obtain
\[
\ll \frac{1}{\sqrt{n}} \int_{-\vartheta_{h,k}'}^{\vartheta_{h,k}''}\sum_{\substack{m_1 \geq 1 \\ m_2 \geq m_p}} |p(m_1)||p(m_2)|e^{\frac{\pi}{12}A(p,m_1,m_2)} d\Phi.
\]
We can bound the integrand from above by $O(1)$ as the series converges absolutely and is independent of $n$ or $k$. If we now plugin the bounds $|\vartheta_{h,k}'|,|\vartheta_{h,k}''| \ll \frac{1}{kN}$, as well as $N=\sqrt{n}$, we find
\[
    \frac{1}{\sqrt{n}}\int_{-\vartheta_{h,k}'}^{\vartheta_{h,k}''} \sum_{m_1,m_2 \geq 0} \zeta^{m_2} p(m_1)p(m_2) e^{\frac{2\pi ih'}{k}(m_1+p^{-1}m_2)} e^{\frac{\pi}{12k}(A(p,m_1,m_2)z^{-1}+(24n-p-1)z)} d\Phi \ll \frac{1}{\sqrt{n}}\cdot \frac{1}{kN} = \frac{1}{n}. 
\]
Now we can take the sum over all Farey arcs and find that the overall contribution is bounded by
\[
\ll \frac{1}{n} \sum_{\substack{0 \leq h < k \leq N \\ \gcd(h,k)=1 \\
    p \nmid k}} e^{-\frac{2\pi ihn}{k} }\sqrt{p} \omega_{h,k}\omega_{hp,k} = \frac{1}{n} \sum_{1 \leq k \leq N} \sum_{\substack{0\leq h < k \\ \gcd(h,k)=1}} e^{-\frac{2\pi ihn}{k} }\omega_{h,k}\omega_{hp,k} \ll \frac{kN}{n} \ll 1.
\]

For the remaining finitely many $m_2$ such that $A(p,0,m_2)>0$ we will use a general result taken from \cite{Bridges-Bringmann} which was shown first in a calculation by Lehner \cite[pp. 404-405]{Lehner}.

\begin{lemma}\label{lemma:integral}
    Suppose that $k\in \N, \nu \in \R$ and $\vartheta_1,\vartheta_2,A,B > 0$ satisfy $k \ll \sqrt{n}, A \ll \frac{n}{k}, B \ll \frac{1}{k},$ and $k\vartheta_1, k\vartheta_2 \asymp \frac{1}{\sqrt{n}}$. Then, we have
    \[
    \int_{\frac{k}{n}-ik\vartheta_1}^{\frac{k}{n}+ik\vartheta_2} Z^{-\nu} e^{AZ+BZ^{-1}} dZ = 2\pi i\left(\frac{A}{B}\right)^{\frac{\nu - 1}{2}} I_{\nu - 1}\left(2\sqrt{AB}\right) + \begin{cases}
        O\left(n^{\nu - \frac{1}{2}}\right), & \nu \geq 0,\\
        \noalign{\vskip5pt}
        O\left(n^{\frac{\nu-1}{2}}\right), & \nu < 0.
    \end{cases} 
    \]
\end{lemma}
Here $I_\nu(z)$ denotes the Bessel function of order $\nu$, see for example \cite{Watson}. We can apply this result almost immediately to the integrals
\[
    \int_{-\vartheta_{h,k}'}^{\vartheta_{h,k}''} z  
    e^{\frac{\pi}{12k}(A(p,0,m_2)z^{-1}+(24n-p-1)z)} d\Phi,
\]
when $A(p,0,m_2) > 0$. To bring it into the shape required we will have to perform the change of variables from $\Phi$ to $z$ completely by taking into account the differential $dz = -ik d\Phi$ and changing the boundaries of the integral. This yields
\[
    =  \frac{1}{ik}\int_{\frac{k}{n}- ik\vartheta_{h,k}''}^{\frac{k}{n}+ik\vartheta_{h,k}'} z
    e^{\frac{\pi}{12k}(A(p,0,m_2)z^{-1}+(24n-p-1)z)} d\Phi.
\]
Now note that since we chose $N=\sqrt{n}$ and $n > \frac{p+1}{24}$ (which guarantees $A>0$ below) in the beginning, we can apply the Lemma with
\[\nu = -1,\quad A=\frac{\pi}{12k}(24n-p-1),\quad B=\frac{\pi}{12k}A(p,0,m_2), \quad \vartheta_1 = \vartheta_{h,k}'',\quad \vartheta_2 = \vartheta_{h,k}'
\]
to obtain (note that $I_{-2}(z)=I_2(z)$ for all $z \in \C$)
\begin{multline*}
\int_{-\vartheta_{h,k}'}^{\vartheta_{h,k}''} z  e^{\frac{\pi}{12k}(A(p,0,m_2)z^{-1}+(24n-p-1)z)} d\Phi \\
= \frac{2\pi}{k}\left(\frac{A(p,0,m_2)}{24n-p-1} \right)I_2\left(\frac{\pi}{6k}\sqrt{A(p,0,m_2)(24n-p-1)}\right) + O\left(\frac{1}{n}\right).
\end{multline*}
We now define the Kloosterman sum
\[K_1(k,p,m_1,m_2,n) \coloneqq \sum_{\substack{0 \leq h < k \\ \gcd(h,k)=1}}\omega_{h,k}\omega_{hp,k}e^{\frac{2\pi i}{k}(h'(m_1+p^{-1}m_2)-hn)} \zeta^{m_2}.
\]
We can thus write \eqref{eqn:contribution-p-coprime-to-k} as
\begin{multline*}
    \sum_{\substack{0 \leq h < k \leq N \\ \gcd(h,k)=1 \\
    p \nmid k}} e^{-\frac{2\pi ihn}{k} }\sqrt{p} \omega_{h,k}\omega_{hp,k} \int_{-\vartheta_{h,k}'}^{\vartheta_{h,k}''} z e^{\frac{\pi}{12k}\left(\left(1+\frac{1}{p}\right)\frac{1}{z}-(1+p)z\right)}e^{\frac{2\pi nz}{k}}P(q_1)P\left(\zeta q_1^{\frac{1}{p}}\right) d\Phi\\ = \sqrt{p}\sum_{\substack{1\leq k \leq N \\ p \nmid k}} \sum_{m_2 = 0}^{m_p} K_1(k,p,0,m_2,n) \frac{2\pi}{k}\left(\frac{A(p,0,m_2)}{24n-p-1} \right)I_2\left(\frac{\pi}{6k}\sqrt{A(p,0,m_2)(24n-p-1)}\right) + O\left(1\right)
\end{multline*}
We now immediately see that the term that grows the fastest comes from $k=1$ and $A(p,0,0)$. We do want to estimate the remainder of the sum against the second fastest growing summand. We shall not worry about the exact shape, and thus introduce the notation
\[
M_1(p) \coloneqq \max_{\substack{k \geq 2, m_2 \geq 0 \\ p \nmid k}} \frac{\sqrt{A(p,0,m_2)}}{k}.
\]
It is clear, that $M_1(p) < \sqrt{A(p,0,0)}$. Hence, we can finally write
\begin{multline*}
    \sum_{\substack{0 \leq h < k \leq N \\ \gcd(h,k)=1 \\
    p \nmid k}} e^{-\frac{2\pi ihn}{k} }\sqrt{p} \omega_{h,k}\omega_{hp,k} \int_{-\vartheta_{h,k}'}^{\vartheta_{h,k}''} z e^{\frac{\pi}{12k}\left(\left(1+\frac{1}{p}\right)\frac{1}{z}-(1+p)z\right)}e^{\frac{2\pi nz}{k}}P(q_1)P\left(\zeta q_1^{\frac{1}{p}}\right) d\Phi\\
    = K\frac{2\pi\sqrt{p} A(p,0,0)}{24n-p-1} I_2\left(\frac{\pi}{6}\sqrt{A(p,0,0)(24n-p-1)}\right) + O\left(\frac{N}{24n-p-1} I_2\left(\frac{\pi}{6}M_1(p)\sqrt{24n-p-1}\right)\right),
\end{multline*}
where we abbreviated $K=K(0,p,0,0,n)$. We now note that $K(0,p,0,0,n)=1$, since as stated before $\omega_{0,1} = 1$. Additionaly, we recall once again that $N=\sqrt{n}$ and compute $A(p,0,0) = 1+p^{-1}$. Thus, the last line from above equals
\[
\frac{2\pi\sqrt{p}}{24n-p-1}I_2\left(\frac{\pi}{6}\sqrt{(1+p^{-1})(24n-p-1)}\right) + O\left(\frac{\sqrt{n}}{24n-p-1}I_2\left(\frac{\pi}{6}M_1(p)\sqrt{24n-p-1}\right)\right).
\]
This ends our computation for the case $p \nmid k$. We mention that the case $p \mid k$ is quite similar and we strongly motivate why this contribution will be negligible. The transformation law in the case $p \mid k$ requires us to consider
\[
 \sum_{\substack{0 \leq h < k \leq N \\ \gcd(h,k)=1 \\
    p \mid k}} e^{-\frac{2\pi ihn}{k} }\sqrt{p} \omega_{h,k}\omega_{h,\frac{k}{p}} \int_{-\vartheta_{h,k}'}^{\vartheta_{h,k}''} z e^{\frac{\pi}{12k}\left((1+p)\frac{1}{z}-(1+p)z\right)}e^{\frac{2\pi nz}{k}}P(q_1)P\left(\zeta q_1^{\frac{1}{p}}\right) d\Phi.
\]
The first summand and largest summand that appears is $k=p$ and if we develop this term in the same way as above we see that it is of highest order
\[
\ll I_2\left(\frac{\pi}{6p}\sqrt{(1+p)(24n-p-1)}\right).
\]
Comapring the argument of the two Bessel functions we see that this term is negligble, since
\[
\frac{1}{p}\sqrt{1+p}=\sqrt{\frac{1+p}{p^2}} = \sqrt{\frac{1+p^{-1}}{p}} < \sqrt{1+p^{-1}}.
\]
and is thus exponentially smaller than the main term from the sum $p\nmid k$. We leave out the remaining details. In summary we find that for all primes $p$ we have as $n \rightarrow \infty$,
\[
a_p(n) \sim 2\pi\sqrt{p}\left(\frac{1+p^{-1}}{24n-p-1}\right)I_2\left(\frac{\pi}{6}\sqrt{(1+p^{-1})(24n-p-1)}\right).
\]
Related but weaker asymptotics were observed in \cite{Kotesovec}. Guadalupe established the above asymptotics in \cite{Guadalupe} for $p<24$ and showed his conjecture for those primes using a computation which holds much more general for any positive integer instead of a prime $p$ and will be reproduced here (from now on until \eqref{eqn:asymptotics-log(a_p(n))}) for completeness.

In the next step we are going to plug in the asymptotic expansion of the Bessel function $I_2(x)$ for large positive $x$. The asymptotics for the modified Bessel functions are well-known (see \cite{Watson} for an encyclopedic reference on Bessel functions, or \cite{Banerjee} for a reference specialized on the asymptotics for large arguments) and we have as $x \rightarrow \infty$,
\[
I_2(x) = \frac{e^x}{\sqrt{2\pi x}}\left(1-\frac{15}{8x} + O\left(\frac{1}{x^2}\right)\right).
\]
Combining these asymptotics we infer
\[
a_p(n) \sim 2\sqrt{3p} \cdot \frac{(1+p^{-1})^{\frac{3}{4}}}{(24n-p-1)^{\frac{5}{4}}}\cdot e^{\frac{\pi}{6}\sqrt{(1+p^{-1})(24n-p-1)}} \left(1 - \frac{45}{8\pi\sqrt{6n(1+p^{-1})}} + O\left(\frac{1}{n}\right)\right).
\]
Furthermore, we expand the remaining terms as follows
\begin{align*}
\frac{1}{(24n-p-1)^{\frac{5}{4}}} &= (24n)^{-\frac{5}{4}}\left(1 + O \left(\frac{1}{n}\right)\right),\\
e^{\frac{\pi}{6}\sqrt{(1+p^{-1})(24n-p-1)}} &=e^{\pi\sqrt{\frac{2n(1+p^{-1})}{3}}}\left(1 - \frac{p^{-\frac{1}{2}} (1+p)^{\frac{3}{2}}\pi}{24\sqrt{6n}} + O\left(\frac{1}{n}\right)\right).
\end{align*}
Therefore, as $n\rightarrow \infty$ the following asymptotics hold
\begin{multline}\label{eqn:asymptotics-a_p(n)}
a_p(n) \sim 2\sqrt{3p} \frac{(1+p^{-1})^{\frac{3}{4}}}{(24n)^{\frac{5}{4}}} \left(1+O\left(\frac{1}{n}\right)\right)e^{\pi\sqrt{\frac{2n(1+p^{-1})}{3}}}\\
\cdot\left(1 - \frac{p^{-\frac{1}{2}} (1+p)^{\frac{3}{2}}\pi}{24\sqrt{6n}} + O\left(\frac{1}{n}\right)\right) \left(1 - \frac{45}{8\pi\sqrt{6n(1+p^{-1})}} + O\left(\frac{1}{n}\right)\right).
\end{multline}
We now take the logarithm on both sides and keep the asymptotic equality. We briefly justify this step. Suppose that $f$ and $g$ are both positive functions bounded from below by some constant $c>1$ and which satisfy $f\sim g$. Then, we have also $\log (f) \sim \log (g)$ since
\[
\lim_{n\rightarrow \infty} \frac{\log (f(n))}{\log (g(n))} = \lim_{n\rightarrow \infty} \frac{\log\left(\frac{f(n)}{g(n)}g(n)\right)}{\log(g(n))} = 1 + \lim_{n\rightarrow \infty} \frac{\log\left(\frac{f(n)}{g(n)}\right)}{\log(g(n))} = 1.
\]
Applying this to \eqref{eqn:asymptotics-a_p(n)} and using the power series expansion of $\log(1+x)$ for $|x|<1$ we do find, as desired,
\begin{equation}\label{eqn:asymptotics-log(a_p(n))}
\log \left(a_p(n)\right) \sim \pi \sqrt{\frac{2n(1+p^{-1})}{3}}-\frac{5}{4}\log \left(n\right) + \log \left(\frac{2\sqrt{3p}(1+p^{-1})^{\frac{3}{4}}}{24^{\frac{5}{4}}}\right) - \frac{c_p}{24\sqrt{6n}}.
\end{equation}
This finishes the proof of \cref{thm:main-theorem} for $r=p$ prime. It is not hard to generalize the proof to general positive $r$. The first step is to note that one will have to consider as many transformation laws as $r$ has divisors. More precisely, we do have (see \cite{Bringmann-Mahlburg}) the following transformation laws
\[
P(q^r) = \omega_{h\frac{r}{\gcd(r,k)},\frac{k}{\gcd(r,k)}}\left(\frac{r}{\gcd(r,k)}z\right)^{\frac{1}{2}}e^{\frac{\pi\gcd(r,k)}{12k}\left(\frac{\gcd(r,k)}{rz}-\frac{r}{\gcd(r,k)}z\right)}P\left(\zeta q_1^{\frac{\gcd(r,k)^2}{r}}\right),
\]
where as before $\zeta=\zeta(r,k,h)$ is a certain root of unity whose exact definition does not matter for our purposes. Observe that for any $k$ with $\gcd(r,k)=1$ (in particular $k=1$) we recover the same transformation law as in the case $r=p$ prime with $k \nmid p$
\[
   P(q^r) = \omega_{hr,k} \sqrt{rz} e^{\frac{\pi}{12k}((rz)^{-1}-rz)}P\left(\zeta q_1^{\frac{1}{r}}\right),
\]
Thus, the same proof as before gives the main term and to show the remaining asymptotics in \cref{thm:main-theorem} it suffices to show that the terms coming from the other transformation laws are exponentially smaller (similar to the case $p \mid k$ from before). Denote $d\coloneqq d(r,k)\coloneqq\gcd(r,k)$ and suppose $d>1$. For each transformation law the case $k=d$ gives the main term of the sum
\[
\sum_{\substack{0 \leq h < k \leq N \\ \gcd(h,k)=1 \\
    \gcd(r,k) = d}} e^{-\frac{2\pi ihn}{k} }\sqrt{r} \omega_{h,k}\omega_{h\frac{r}{d},\frac{k}{d}} \int_{-\vartheta_{h,k}'}^{\vartheta_{h,k}''} z e^{\frac{\pi}{12k}\left(\left(1+\frac{d^2}{r}\right)\frac{1}{z}-(1+r)z\right)}e^{\frac{2\pi nz}{k}}P(q_1)P\left(\zeta q_1^{\frac{d^2}{r}}\right) d\Phi.
\]
and the contribution is by \cref{lemma:integral} bounded (as in the case $r=p$) by
\[
\ll I_2\left(\frac{\pi}{6d}\sqrt{\left(1+ \frac{d^2}{r}\right)(24n-r-1)}\right)
\]
The critical part of the argument inside the Bessel function is bounded by
\[
\frac{1}{d}\sqrt{1+\frac{d^2}{r}}=\sqrt{\frac{1}{d^2}\left(1+\frac{d^2}{r}\right)} = \sqrt{\frac{1}{d^2}+\frac{1}{r}} < \sqrt{1+\frac{1}{r}}.
\]
The right hand side is the quantity from the case $d=1$. Thus, we have shown that all main terms coming from the $d>1$ parts of the sum are exponentially smaller than main term coming from $d=1$. Therefore, we have
\[
a_r(n) \sim 2\pi\sqrt{r}\left(\frac{1+r^{-1}}{24n-r-1}\right)I_2\left(\frac{\pi}{6}\sqrt{(1+r^{-1})(24n-r-1)}\right).
\]
The rest of the proof follows mutas mutandis to the case $r=p$ prime. \hfill \qedsymbol

\bibliographystyle{alpha}

\end{document}